\documentclass[parskip=full, twoside=false]{scrartcl} 
\usepackage[utf8]{inputenc} 
\usepackage[T1]{fontenc} 
\usepackage[english]{babel}
\usepackage[autostyle=true]{csquotes}
\usepackage{lmodern} 
\usepackage{cite}
\usepackage[pdfauthor={Jonathan Glöckle}, pdftitle={An Enlargeability Obstruction for Spacetimes with both Big Bang and Big Crunch}]{hyperref} 
\usepackage{amsmath, amsthm, amssymb, bm, bbm, mathrsfs, mathtools}
\usepackage{extarrows}
\usepackage{tabulary,tabularx}
\usepackage{graphicx}
\usepackage{tikz-cd}
\usepackage[capitalise]{cleveref}
\usepackage{bbold}

\newcommand{\arxiv}[1]{\href{https://arxiv.org/abs/#1}{ArXiV:#1}}

\theoremstyle{plain} 
\newtheorem{Satz}{Theorem}[section] 
\newtheorem{Prop}[Satz]{Proposition} 
\newtheorem{Lem}[Satz]{Lemma}
\newtheorem{Kor}[Satz]{Corollary} 
\newtheorem*{Beh*}{Claim}
\newtheorem*{Satz*}{Main Theorem}
\theoremstyle{definition} 
\newtheorem{Def}[Satz]{Definition} 

\newtheorem{Bem}[Satz]{Remark}

\crefname{Satz}{Theorem}{Theorems}
\crefname{Prop}{Proposition}{Propositions}
\crefname{Lem}{Lemma}{Lemmata}
\crefname{Kor}{Corollary}{Corollaries}
\crefname{Def}{Definition}{Definitions}

\newcommand{\N}{\mathbbm{N}} 
\newcommand{\Z}{\mathbbm{Z}} 
 
\newcommand{\R}{\mathbbm{R}} 
\newcommand{\C}{\mathbbm{C}}

\newcommand{\del}{\partial}

\newcommand{\lto}{\longrightarrow}
\newcommand{\lmapsto}{\longmapsto}

\newcommand{\upd}{\mathrm{d}}

\newcommand{\Cl}{\C\mathrm{l}}
\newcommand{\pt}{\{*\}}

\renewcommand{\epsilon}{\varepsilon}
\renewcommand{\theta}{\vartheta}

\let\div\relax

\DeclareMathOperator{\PSC}{\mathcal{R}^{>}}
\DeclareMathOperator{\Ini}{\mathcal{I}}
\DeclareMathOperator{\DEC}{\mathcal{I}^{>}}

\DeclareMathOperator{\ind}{ind}
\DeclareMathOperator{\idiff}{ind-diff}

\DeclareMathOperator{\ric}{ric} 
\DeclareMathOperator{\div}{div}
\DeclareMathOperator{\scal}{scal}
\DeclareMathOperator{\Ein}{Ein}

\DeclareMathOperator{\tr}{tr}
\DeclareMathOperator{\id}{\mathbb{1}}

\DeclareMathOperator{\supp}{supp}
\DeclareMathOperator{\SO}{SO}
\DeclareMathOperator{\Spin}{Spin}

\DeclareMathOperator{\Ind}{Ind}

\DeclareMathOperator{\A}{\hat A-deg}
\DeclareMathOperator{\KO}{KO}

\begin{document}
\title{An Enlargeability Obstruction for Spacetimes with both Big Bang and Big Crunch}
\author{Jonathan Glöckle}
\maketitle

\begin{abstract}
	Given a spacelike hypersurface $M$ of a time-oriented Lorentzian manifold $(\overline{M},\overline{g})$, the pair $(g,k)$ consisting of the induced Riemannian metric $g$ and the second fundamental form $k$ is known as initial data set.
	In this article, we study the space of all initial data sets $(g,k)$ on a fixed closed manifold $M$ that are subject to a strict version of the dominant energy condition.
	Whereas the pairs of the form $(g,\tau g)$ and $(g,-\tau g)$, for a sufficiently large $\tau > 0$, belong to the same path-component of this space when $M$ admits a positive scalar curvature metric, it was observed in a previous work \cite{gloeckle:p2019} that this is not the case when the existence of a positive scalar curvature metric on $M$ is obstructed by $\alpha(M) \neq 0$.
	In the present article we extend this non-connectedness result to Gromov-Lawson's enlargeability obstruction, which covers many examples, also in dimension $3$.
	In the context of relativity theory, this result may be interpreted as excluding the existence of certain globally hyperbolic spacetimes with both a big bang and a big crunch singularity.
\end{abstract}

\section{Introduction}
The interplay of curvature restrictions and the topology of a manifold has always been a central topic in differential geometry.
Curvature restrictions also have important applications in general relativity.
Namely, the famous singularity theorems of Hawking and Penrose (cf.~\cite[Sec.~8.2]{hawking.ellis:73}) rely on the fact that all matter is assumed to satisfy certain energy conditions, which, by the Einstein field equations, correspond to curvature conditions for the Lorentzian space-time manifold.
For example, the dominant energy condition -- derived from the physical assumption that every observer should experience non-negative mass density -- translates into a “non-negativity” condition for the Einstein curvature $\Ein = \ric - \frac{1}{2} \scal g$.

A usual strategy for studying solutions of the Einstein equations is to look at time slices:
Let $M$ be a spacelike hypersurface of a time-oriented Lorentzian manifold $(\overline{M}, \overline{g})$.
The Lorentzian metric $\overline{g}$ induces on $M$ a Riemannian metric $g$ and a symmetric $2$-tensor $k \in \Gamma(T^*M \otimes T^*M)$, the second fundamental form of $M$ with respect to the future-directed unit normal $e_0$.
The Gauß-Codazzi equations imply that the components $\rho = \Ein(e_0,e_0)$ and $j = \Ein(e_0, -) \in \Omega^1(M)$ of the Einstein curvature may be computed via the so-called \emph{constraint equations}
\begin{equation} \label{eq:CE}
\begin{aligned}
	2\rho &= \scal^{g} + (\tr k)^2 - \|k\|^2 \\
	j &= \div k - \upd \tr k.
\end{aligned} 
\end{equation}
By the Einstein equations, $\rho$ and $j$ are determined by the distribution of matter: $\rho$ is the \emph{energy density} and $j$ gets the interpretation of a \emph{momentum density}.
If $\overline{g}$ satisfies the dominant energy condition, then $\rho \geq \|j\|$ everywhere on $M$.
At least a partial converse of this is also true:
Suppose that a pair $(g,k)$ of metric and symmetric $2$-tensor on $M$ is such that $\rho$ and $j$ -- determined by \eqref{eq:CE} -- satisfy the strict inequality $\rho > \|j\|$ on all of $M$.
Then $M$ embeds as a spacelike hypersurface into some time-oriented Lorentzian manifold $(\overline{M}, \overline{g})$ subject to the dominant energy condition such that $g$ is the induced metric and $k$ is the induced second fundamental form.
Therefore, we say that an \emph{initial data set} $(g, k)$ satisfies the \emph{dominant energy condition} (dec) if $\rho \geq \|j\|$.
This condition plays a major role in the positive mass theorem.

If $M \subseteq \overline{M}$ is totally geodesic, i.\,e.~if $k  \equiv 0$, then the dominant energy condition reduces to the condition that $g$ has non-negative scalar curvature.
Non-negative and even more so positive scalar curvature (psc) has been a vast field of study over the last decades, see e.\,g.~\cite{rosenberg:2007}.
Apart from minimal hypersurface techniques, most of these results were obtained by Dirac operator methods.
For example, it was shown in \cite{botvinnik.ebert.randal-williams:14} and \cite{crowley.steimle.schick:18} that the $C^\infty$-space $\PSC(M)$ of psc metrics on a fixed closed spin manifold $M$ of dimension $n \geq 6$ has infinitely many non-trivial homotopy groups if it is non-empty.

In \cite{gloeckle:p2019}, the author was able to transfer some of these psc results to the dominant energy setting.
Namely, the non-trivial elements in $\pi_k(\PSC(M))$ mentioned above give rise to non-trivial elements in $\pi_{k+1}(\DEC(M))$, where $\DEC(M)$ denotes the space of initial data sets $(g, k)$ on $M$ subject to the strict version $\rho > \|j\|$ of the dominant energy condition.
The proof relies on a suspension construction -- associating to a psc metric a path of initial data sets strictly satisfying dec -- and the computation of a family index of (Clifford-linear) \emph{Dirac-Witten operators}.
These operators are a certain zero-order perturbation of Dirac operators.
Their classical version dates back to Witten's spinorial proof of the positive mass theorem \cite{witten:81} (cf.~\cite{parker.taubes:82} for a more detailed account of the proof).

In a some way, the comparison arguments between $\pi_k(\PSC(M))$ and $\pi_{k+1}(\DEC(M))$ still apply in the case $k = -1$.
Here, statement is that $\DEC(M)$ is not path-connected if $M$ is a closed spin manifold with non-zero $\alpha$-index $\alpha(M)$.
More precisely, the path-component $C^+$ that contains all data sets $(g, \tau g) \in \DEC(M)$ for sufficiently large $\tau > 0$ and the path-component $C^-$ that contains $(g, -\tau g)$ for large $\tau > 0$ are distinct in this case.

This has a direct application to general relativity.
Near a big bang singularity the initial data sets induced on spacelike hypersurfaces are expected to lie in $C^+$, whereas those near a big crunch singularity would belong to $C^-$.
Recall that a globally hyperbolic spacetime $(\overline{M}, \overline{g})$ admits a foliation $\overline{M} \cong M \times \R$ into spacelike (Cauchy) hypersurfaces $M \times \{t\}$ by a famous theorem by Bernal and Sánchez \cite{bernal.sanchez:05}.
Taking this together, we see that $\alpha(M) \neq 0$ is not only an obstruction to psc but also an obstruction to the existence of a globally hyperbolic spacetime $(\overline{M} \cong M \times \R, \overline{g})$ strictly satisfying dec and having both a big bang and a big crunch singularity.
The slight drawback here -- that we require $(\overline{M}, \overline{g})$ to strictly satisfy dec meaning that on every spacelike hypersurface $\rho > \|j\|$ holds -- is overcome in the recent article \cite{ammann.gloeckle:p2021} of Ammann and the author, where sufficient conditions for the removal of the strictness assumption are deduced.
This passage from strict dec to dec is similar to the passage from positive scalar curvature to non-negative scalar curvature performed by Schick and Wraith \cite{schick.wraith:2021}.

Unfortunately, precisely in the physically relevant spatial dimension $n = 3$, we always have $\alpha(M) = 0$ as it is an element of $\KO^{-3}(\pt) = 0$.
It is the purpose of the present article to extend the result of \cite{gloeckle:p2019} to another (index theoretic) obstruction to positive scalar curvature that is also of use in the $3$-dimensional case: enlargeability.
The concept of enlargeability was first introduced by Gromov and Lawson \cite{gromov.lawson:1980b}.
There exist various versions, the most basic ones being the following.
A closed Riemannian $n$-manifold $(M, g)$ is called \emph{enlargeable} if for every $R > 0$ there exists a Riemannian covering of $(M, g)$ that is “large” in the following sense:
There exists a distance-non-increasing map to a standard sphere $S^n$ of radius $R$ that is of non-zero degree.
Despite being geometrically defined, the notion turns out to be homotopy invariant (in particular, independent of $g$), and thus provides a link between geometry and topology.
Moreover, the class of enlargeable manifolds is rich.
It contains all manifolds that carry a metric of non-positive sectional curvature, especially tori, and further examples may be constructed through products and direct sums.
In the aforementioned article \cite{gromov.lawson:1980b} Gromov and Lawson show that spin manifolds do not carry a psc metric when they are \emph{compactly enlargeable}, i.\,e.~the coverings in the definition of enlargeability may be chosen to be finite-sheeted.
Later they extended this both to (not necessarily compactly) enlargeable spin manifolds and, in dimensions $n \leq 7$, to compactly enlargeable (not necessarily spin) manifolds (cf.~\cite{gromov.lawson:1983}).
In the present article we prove an initial data version of Gromov and Lawson's result:
\begin{Satz*}[\Cref{Thm:EnlObst}]
If $M$ is an enlargeable spin manifold, then $\DEC(M)$ is not path-connected.
More precisely, the path-components $C^+$ and $C^-$ are distinct.
\end{Satz*}
In fact, we prove this for the slighly more general notion of $\hat{A}$-area-enlargeability, cf.~\cref{Def:Enl}.
Similar to what has been stated above, this has the following consequence to general relativity:
If $M$ is an enlargeable spin manifold and $(\overline{M} \cong M \times \R,\overline{g})$ a globally hyperbolic spacetime that strictly satisfies dec in the sense that $\rho > \|j\|$ on every spacelike hypersurface, then $(\overline{M}, \overline{g})$ cannot have both a big bang and a big crunch singularity.
In particular, this applies when $M$ is a quotient of a ($3$-dimesional) torus, which is still considered to be within the range of physical observations of our universe, cf.~\cite{galloway.khuri.woolgar:p2020}.
Probably, the strictness condition can be removed under certain additional assumptions with techniques as in \cite{ammann.gloeckle:p2021} and \cite{schick.wraith:2021}, but this is beyond the scope of this article.

Let us briefly discuss the idea of the proof in the case where $M$ is even-dimensional and compactly enlargeable.
From the (compact) enlargeability Gromov and Lawson construct a sequence of (finite-sheeted) coverings $M_i \to M$ and complex vector bundles $E_i \to M_i$ with hermitian metric and metric connection such that the curvatures $R^{E_i} \lto 0$ for $i \lto \infty$ and $\hat{A}(M_i, E_i) \neq 0$ for all $i \in \N$.
Thereby $\hat{A}(M_i, E_i) = \int_{M_i} \hat{A}(TM_i) \wedge \mathrm{ch}(E_i)$ is the $E_i$-twisted $\hat{A}$-genus, which is equal to the index of the $E_i$-twisted Dirac operator $D^{E_i}$ by the Atiyah-Singer index theorem.
Note that the existence of such a sequence of coverings and bundles is also the starting point of Hanke and Schick's proof that enlargeability implies non-triviality of the Rosenberg index $\alpha_{\mathrm{max}}^\R(M) \in \KO^{-n}(C^*_{\mathrm{max},\R}\pi_1(M))$ \cite{hanke.schick:2006,hanke.schick:2007}.
In our initial data setting, we need a suitable analog of $D^{E_i}$.
This is provided by the twisted Dirac-Witten operator $\overline{D}^{E_i} = D^{E_i} -\frac{1}{2} \tr(k) e_0 \cdot$ defined on a certain $E_i$-twisted spinor bundle. 
This operator is associated to an initial data set $(g, k)$ on $M_i$ (in our case it is pulled back from $M$) and satisfies the Schrödinger-Lichnerowicz type formula
\begin{align*}
	\left(\overline{D}^{E_i}\right)^2 \psi = \overline{\nabla}^* \overline{\nabla} \psi + \frac{1}{2} (\rho - e_0 \cdot j^\sharp \cdot) \psi + \mathcal{R}^{E_i} \psi.
\end{align*} 
Hence, if the strict dominant energy condition $\rho > \|j\|$ holds, then $\overline{D}^{E_i}$ is invertible for large enough $i \in \N$.
In particular, given a metric $g$ and a sufficiently large $\tau > 0$ the \emph{twisted index difference} $\idiff^{E_i}((g, -\tau g), (g, \tau g))$, a spectral-flow-like invariant associated to a path $\gamma$ of initial data sets from $(g, -\tau g)$ to $(g, \tau g)$, vanishes for large $i \in \N$ if the path $\gamma$ may be chosen to lie within $\DEC(M)$.
On the other hand, an index theorem shows that $\idiff^{E_i}((g, -\tau g), (g, \tau g)) = \hat{A}(M_i, E_i) \neq 0$ for all $i \in \N$.
Thus it is not possible to connect $(g, -\tau g)$ and $(g, \tau g)$ by a path in $\DEC(M)$.

The article is structured as follows:
We start off by describing how Dirac-Witten operators arise in a Lorentzian setup.
In the second section, we construct the twisted index difference of Dirac-Witten operators using a framework laid out by Ebert \cite{ebert:p2018}.
This is a bit technical, as we also deal with non-compact manifolds so that we can later apply it to infinite covers $M_i \to M$ as well.
In fact, we construct a \emph{relative twisted index difference}, relative meaning that it depends on the “difference” of two twist bundles $E_i^{(0)}$ and $E_i^{(1)}$ that coincide outside a compact set.
The third section is devoted to the proof of the (relative) index theorem $\idiff^{E_i}((g, -\tau g), (g, \tau g)) = \hat{A}(M_i, E_i)$.
In the last section we put the arguments together to prove the main theorem.

\subsection*{Acknowledgements}
I would like express my gratitude to Bernd Ammann, who came up with the idea for this project, for his ongoing support and encouragement. During its execution, I was supported by the SFB 1085 \enquote{Higher Invariants} funded by the DFG.

\section{Lorentzian manifolds, initial data sets and Dirac-Witten operators}
Within this section, we want to recall several notions from Lorentzian geometry and thereby fix certain notations.
Although the rest of the article could be understood to a very large extent without knowledge of the Lorentzian setup, the objects that we will be considering appear more naturally in this context, thus providing a better understanding.

Let us consider a Lorentzian manifold $(\overline{M}, \overline{g)}$ of dimension $n+1$.
We use the signature convention such that a generalized orthonormal basis $e_0, e_1, \ldots, e_n$ satisfies $\overline{g}(e_0,e_0) = -1$ and $\overline{g}(e_i,e_i) = 1$ for $1 \leq i \leq n$.
In particular, the induced metric $g$ on a spacelike hypersurface $M \subseteq \overline{M}$ will be Riemannian.
In this context, we will denote by $e_0$ a (generalized) unit normal on $M$.
Usually, we will assume that $\overline{M}$ is time-oriented and then we agree that $e_0$ is future-pointing.
Apart from the induced metric $M$ will carry an induced second fundamental form $k$ with respect to $e_0$ that we define by
\begin{align*}
	\overline{\nabla}_X Y- \nabla_X Y &= k(X,Y) e_0
\end{align*}
for all vectors fields $X, Y \in \Gamma(TM)$.
Thereby, $\overline{\nabla}$ denotes the Levi-Civita connection of $(\overline M, \overline g)$ and $\nabla$ is the one of the hypersurface $(M, g)$.
Pairs $(g, k)$ consisting of a Riemannian metric $g$ and a symmetric $2$-tensor $k$ will be called \emph{initial data set}.
Hence the above procedure provides an induced initial data set $(g, k)$ on a spacelike hypersurface of a time-oriented Lorentzian manifold.

The \emph{dominant energy condition} is the condition that for all future-causal vectors $V, W$ of $(\overline M, \overline{g})$, the Einstein tensor $\Ein = \ric - \frac{1}{2}\scal \overline g$ satisfies $\Ein(V,W) \geq 0$.
Equivalently, for every future-causal vector $V$ the metric dual of $\Ein(V,-)$ is required to be past-causal.
Applying this to the future-pointing unit normal $e_0$ of a spacelike hypersurface $M$, we get $\rho \geq \|j\|$ for $\rho = \Ein(e_0, e_0)$ and $j = \Ein(e_0, -) \in \Omega^1(M)$.
As explained in the introduction, $\rho$ and $j$ are completely determined by the induced initial data set $(g, k)$ on $M$ via \eqref{eq:CE}.

\begin{Def} \label{Def:DEC}
An initial data set $(g,k)$ on a manifold $M$ is said to satisfy the \emph{dominant energy condition} if $\rho \geq \|j\|$ for $\rho$ and $j$ defined by \eqref{eq:CE}.
It satisfies the \emph{strict dominant energy condition} if $\rho > \|j\|$.
When $M$ is compact, we denote by $\Ini(M)$ the space of all initial data sets $(g, k)$ on $M$, equipped with the $C^\infty$-topology of uniform convergence, and by $\DEC(M)$ the subspace of those initial data sets strictly satisfying the dominant energy condition.
\end{Def}
It should be noted that when $M$ is compact of dimension $n \geq 2$ and $g$ is any metric on $M$, then $(g, \tau g)$ and $(g, -\tau g)$ satisfy the strict dominant energy condition once the function $\tau > 0$ is large enough.
In fact, in this case $\rho > \|j\|$ reduces to
\begin{align} \label{eq:DECsp}
	\frac{n(n-1)}{2} \tau^2 > -\frac{1}{2} \scal^g + (n-1) \|\upd \tau\|,
\end{align}
which can always be achieved by adding a constant to $\tau$ that continuously depends on $g$ (in $C^\infty$-topology).
As, moreover, for fixed $g$ the convex combination between some $\tau > 0$ and the constant $\max(\tau)$ keeps satisfying \eqref{eq:DECsp}, all pairs of the form $(g, \tau g)$ for $\tau > 0$ belong to the same path-component $C^+$ of $\DEC(M)$, likewise $(g, -\tau g)$ all belong to the same path-component $C^-$.
When there exists a positive scalar curvature metric $g$ on $M$, then the path $[-1,1] \to \DEC(M),\,t \mapsto (g, t g)$ shows that $C^+ = C^-$.
For future reference, we also note that the question of whether $C^+ = C^-$ satisfies a certain stability property:

\begin{Lem} \label{Lem:Stab}
	Let $M$ be a compact manifold such that the path-components $C^+$ and $C^-$ of $\DEC(M)$ coincide.
	Then for any compact manifold $N$ the path-components $C^+$ and $C^-$ of $\DEC(M \times N)$ agree as well.
	\begin{proof}
		Let us first consider pairs of product form, i.~e.\  $(g_M + g_N, k_M + k_N)\in \Ini(M \times N)$ for $(g_M,k_M) \in \Ini(M)$ and $(g_N,k_N) \in \Ini(N)$. For these we compute
		\begin{align*}
		2\rho &= \scal^{g_M} + \scal^{g_N} + (\tr^{g_M}(k_M) + \tr^{g_N}(k_N))^2 - \|k_M\|_{g_M}^2 - \|k_N\|_{g_N}^2 \\
		&= 2\rho_M + 2\rho_N + 2\tr^{g_M}(k_M)\tr^{g_N}(k_N) \\
		\intertext{and}
		j &= \div^{g_M}(k_M) + \div^{g_N}(k_N) - d(\tr^{g_M}(k_M)) - d(\tr^{g_N}(k_N)) \\
		&= j_M + j_N.
		\end{align*}
		Hence, we obtain
		\begin{align} \nonumber
		\|j\|^2 &= \|j_M\|_{g_M}^2 + \|j_N\|_{g_N}^2 \leq (\|j_M\|_{g_M} + \|j_N\|_{g_N})^2
		\intertext{and thus} \label{eq:prodpairs}
		\rho - \|j\| &\geq (\rho_M - \|j_M\|_{g_M}) + (\rho_N - \|j_N\|_{g_N}) + \tr^{g_M}(k_M)\tr^{g_N}(k_N).
		\end{align}
		
		Now, by assumption, for sufficiently large $\tau > 0$ the pairs $(g_M, -\tau g_M)$ and $(g_M, \tau g_M)$ can be connected by a path $t \mapsto \gamma_M(t)=(g_M (t),k_M (t))$ in $\DEC(M)$. By the discussion following \Cref{Def:DEC}, we may assume that $\tau$ is a constant. As the interval $[-1,1]$ is compact, $\rho_M - \|j_M\|_{g_M}$ attains a positive minimum. Replacing $\gamma_M$ by $\widetilde{\gamma_M} = (C^{-2} g_M, C^{-1} k_M)$ for some suitably chosen $C > 0$, we may assume that this minimum is larger than $-\frac12 \min_{p \in N} \scal^{g_N}(p)$. This is due to the fact, that in this rescaling $\widetilde{\rho_M} = C^2 \rho_M$ and $\left\|\widetilde{j_M}\right\|_{\widetilde{g_M}} = C^2 \|j_M\|_{g_M}$, indicating the rescaled quantities by $\widetilde{\cdot}$.
		
		The required path from $(g_M + g_N, -\tau (g_M + g_N))$ to $(g_M + g_N, \tau (g_M + g_N))$ can be easily pieced together from three segments: 
		\begin{align*}
		[-1,1] &\lto \DEC(M \times N) \\
		t &\lmapsto
		\begin{cases}
		(g_M + g_N, -\tau g_M + (2t+1) \tau  g_N) & t \in [-1,-\frac12] \\
		(g_M(2t) + g_N, k_M(2t)) &t \in [-\frac12,\frac12] \\			
		(g_M + g_N, \tau g_M + (2t-1) \tau g_N) & t \in [\frac12,1].
		\end{cases}
		\end{align*}
		In the first section, both $\tr^{g_M}(k_M)$ and $\tr^{g_N}(k_N)$ are non-positive, so its product is non-negative. Furthermore, for the pair $(g_N, (2t+1)g_N)$ with $t \in [-1, -\frac12]$ we have
		\begin{align*}
		\rho_N - \|j_N\|_{g_N} = \frac{1}{2} (\scal^{g_N} + \dim(N)(\dim(N)-1) (2t+1)^2 \tau^2) \geq \frac12 \min_{p \in N} \scal^{g_N}(p).
		\end{align*}
		By choice of the rescaling, we have that $(g_M, -\tau g_M) = \gamma_M(-1)$ satisfies $\rho_M - \|j_M\|_{g_M} > - \frac12 \min_{p \in N} \scal^{g_N}(p)$, and so \eqref{eq:prodpairs} shows that the first section lies in $\DEC(M \times N)$. The same argument applies for the last section. In the middle section, the last term in \eqref{eq:prodpairs} is zero, $\rho_N - \|j_N\|_{g_N} \geq \frac12 \min_{p \in N} \scal^{g_N}(p)$, and so by our choice of rescaling, the pair is in $\DEC(M \times N)$ for all $t \in [-\frac{1}{2}, \frac12]$.
	\end{proof}
\end{Lem}

We conclude this section by a brief discussion of the Dirac-Witten operator.
Assume that $(\overline M, \overline g)$ is a space- and time-oriented Lorentzian spin manifold.
Let $\Sigma \overline M \to \overline M$ be the classical spinor bundle of $(\overline M, \overline g)$, i.\,e.~the spinor bundle associated to an irreducible representation of $\Cl_{n,1}$.
A short summary of spin geometry in the semi-Riemannian setting can be found in \cite{baer.gauduchon.moroianu:2005}.
Restricting this bundle to the spacelike hypersurface $M$ yields the induced \emph{hypersurface spinor bundle} $\Sigma \overline{M}_{|M} \to M$.
The Levi-Civita connection of $(\overline{M}, \overline{g})$ induces a connection $\overline{\nabla}$ on $\Sigma \overline{M}_{|M}$, and the associated Dirac type operator $\overline{D} \psi = \sum_{i = 1}^n e_i \cdot \overline{\nabla}_{e_i} \psi$, where $(e_1, \ldots, e_n)$ is an local orthonormal frame of $TM$, is known as \emph{Dirac-Witten operator}.
It was first observed by Witten \cite{witten:81} that it is linked to the dominant energy condition by a Schrödinger-Lichnerowicz type formula:
\begin{align} \label{eq:SL}
	\overline{D}^2 \psi = \overline{\nabla}^* \overline{\nabla} \psi + \frac{1}{2}(\rho - e_0 \cdot j^\sharp \cdot) \psi
\end{align}
Here, $\psi$ is any smooth section of $\Sigma \overline{M}$ and $-^*$ denotes the formal adjoint.
The formula shows that the dominant energy condition implies that $\overline{D}$ is positive and hence invertible.

It now turns out, that $\overline{D}$ does not depend on the whole Lorentzian manifold $(\overline{M}, \overline{g})$, but only on induced initial data set $(g,k)$ on $M$.
In fact, denoting by $e_0$ the future-pointing unit normal as above, we obtain an induced $\Spin(n)$-principal bundle on $M$ by pulling back the $\Spin_0(n,1)$-principal bundle $P_{\Spin_0(n,1)} \overline{M}_{|M}$ along the inclusion of frame bundles
\begin{align*}
	P_{\SO(n)}M &\lto P_{\SO_0(n,1)} \overline{M}_{|M} \\
	(e_1, \ldots, e_n) &\lmapsto (e_0, e_1, \ldots, e_n).
\end{align*}
When $n=2m$ is even, there are two irreducible (ungraded) $\Cl_{n,1}$-representations and either restricts to the unique irreducible (ungraded) representation of $\Cl_n$ along the inclusion\footnote{Some authors identify $\Cl_n$ with the subalgebra $\Cl_{n,1}^0$ instead, for instance this is used in \cite{baer.gauduchon.moroianu:2005}.}
$\Cl_n \hookrightarrow \Cl_{n,1}$.
On associated bundles, this yields an isomorphism $\Sigma \overline{M}_{|M} \cong \Sigma M$, where $\Sigma M \to M$ is the classical spinor bundle on $M$.
The difference between the two representations results in the fact that in one case multiplication by $e_0$ is given by $\omega=i^{m} e_1 \cdots e_n$, whereas in the other case it is given by $-\omega$.

Apart from $\overline{\nabla}$, there is another canonical connection on $\Sigma \overline{M}_{|M} \cong \Sigma M$: the one induced by the Levi-Civita connection of $(M,g)$, called $\nabla$.
Those two connections differ by a term depending on the second fundamental form, namely
\begin{align*}
	\overline{\nabla}_X \psi = \nabla_X \psi -  \frac{1}{2} e_0 \cdot k(X,-)^\sharp \cdot \psi
\end{align*}
for $\psi \in \Gamma(\Sigma M)$ and $X \in TM$.
As a consequence, we obtain
\begin{align*}
	\overline{D} \psi = D \psi - \frac{1}{2} \tr(k) e_0 \cdot \psi
\end{align*}
for all $\psi \in \Gamma(\Sigma M)$, where $D$ denotes the Dirac operator on $\Sigma M$. 

Unlike the Dirac operator $D$, which -- for even $n = 2m$ -- is odd with respect to the $\Z/2\Z$-grading defined by the volume element $\omega$ of $\Sigma M$, there is no natural grading which is compatible with the Dirac-Witten operator in this case.
As for index theory, however, such gradings are very useful, we consider a $\Z/2\Z$-graded version instead.
We do so by replacing $\Sigma \overline M$ with $\overline{\Sigma} \overline M$, the bundle associated to the unique irreducible $\Z/2\Z$-graded $\Cl_{n,1}$-representation\footnote{The grading on the Clifford algebra is given by the even-odd-grading.
Note that then a $\Z/2\Z$-graded $\Cl_{n,1}$-representation is nothing else than an ungraded $\Cl_{n,2}$-representation.}.
This representation is obtained by summing the two irreducible ungraded $\Cl_{n,1}$-representations and taking an appropriate grading.
More precisely, starting with the irreducible representation with $i^m e_0 \cdots e_n = \id$, the other one can be obtained by replacing $e_0 \cdot$ with $-e_0 \cdot$, keeping the multiplication by the other basis vectors the same.
Then the grading on the sum can be chosen as
\begin{align} \label{eq:grd}
	\iota &= \begin{pmatrix}
		0 & \omega \\ \omega & 0
	\end{pmatrix}
\end{align}
for $\omega=i^{m} e_1 \cdots e_n$ as above.
As a consequence, there is an isomorphism $\overline{\Sigma} \overline{M}_{|M} \cong \Sigma M \oplus \Sigma M$ such that the grading is given by \eqref{eq:grd} and multiplication by $e_0$ is given by
\begin{align} \label{eq:e0}
	e_0 \cdot \;&= \begin{pmatrix}
		\omega & 0 \\ 0 & -\omega
	\end{pmatrix}.
\end{align}
It should be noted that doubling the spinor bundle creates an additional symmetry:
\begin{align} \label{eq:c1}
	c_1 &= \begin{pmatrix}
		i & 0 \\ 0 & -i
	\end{pmatrix}
\end{align}
defines an odd $\Cl_1$-action on $\overline{\Sigma} \overline{M}_{|M}$, which has the property that it commutes with the Dirac-Witten operator.
For this reason, the Dirac-Witten operator on this bundle is called \emph{$\Cl_1$-linear Dirac-Witten operator}.

\section{The twisted index difference for Dirac-Witten operators}
Let $M$ be a connected spin manifold of even dimension $n=2m$.
We will not assume that $M$ is compact.
The necessity for also considering non-compact manifolds
-- although the main result is only concerned with compact ones --
strives from the fact, that we will later be looking at coverings
and we do not want to assume them to be finite.
The aim of this section is to define a homotopy invariant for a path $\gamma \colon (I, \del I)\to (\Ini(M), \Ini^{> \mathcal{R}^E + c}(M))$, where $I=[-1,1]$ and $c > 0$. This invariant will have the property of being zero if a representative of its homotopy class takes values only in $\Ini^{> \mathcal{R}^E + c}(M)$. The way, it is constucted, is quite similar to the $\overline{\alpha}$-difference in \cite{gloeckle:p2019}, yet it differs in the way that instead of real Clifford-linear spinors, it uses complex spinors with coefficients in a twist bundle. This is needed to make use of enlargeability. 

To be able to define this also when $M$ is non-compact, we need some extra care in the definition of $\Ini(M)$: It will be the space of pairs $(g,k)$, where $g$ is a \emph{complete} Riemannian metric and $k \in \Gamma(T^*M \otimes T^*M)$ is symmetric. It carries the topology of uniform convergence of all derivatives on compact sets. $\Ini^{> \mathcal{R}^E + c}(M)$ denotes the subspace of those pairs satisfying the dominant energy condition in the stricter sense $\frac{1}{2}(\rho - \|j\|) > \|\mathcal{R}^E\| + c$ (cf.\ \eqref{eq:CE}), where $\mathcal{R}^E$ is the curvature endomorphism of a twist bundle $E$.

The twist bundle $E$ arises as a direct sum of two
complex vector bundles $E_0, E_1 \to M$ with hermitian metrics and metric connections, such that outside a compactum $K$ they can be identified by an isometric and connection preserving bundle isomorphism $\Psi \colon E_{0|M \setminus K} \to E_{1|M \setminus K}$. If $M$ is already compact, then, of course, we may take $K = M$ and the compatibility condition becomes void.

The constuction of the desired homotopy invariant begins as follows: For a chosen (complete) Riemannian metric $g$ on $M$, let $\Sigma M$ be the classical complex spinor bundle associated to the (topological) spin structure on $M$. We consider the double spinor bundle $\overline{\Sigma} M \coloneqq \Sigma M \oplus \Sigma M$ with its direct sum hermitian metric. This carries a (self-adjoint) $\Z/2\Z$-grading $\iota$ and a (skew-adjoint) odd $\Cl_1$-action $c_1$ given by \eqref{eq:grd} and \eqref{eq:c1}, respectively. Moreover, $\overline{\Sigma} M$ admits an operator $e_0 \cdot$ defined by \eqref{eq:e0}, which is self-adjoint, odd and commutes with the $\Cl_1$-action. If $D$ is the Dirac operator of the double spinor bundle, then the $\Cl_1$-linear Dirac-Witten operator
\begin{align*}
	\overline{D} = D - \frac{1}{2} \tr^g(k) e_0 \cdot
\end{align*}
is formally self-adjoint, odd and commutes with the $\Cl_1$-action. We mean to associate a suitable (relative) index to it.

We now bring in the twist bundles. From $E_0$ and $E_1$, we form the sum $E = E_0 \oplus E_1$, which we endow with the $\Z/2\Z$-grading
\begin{align*}
	\eta = \begin{pmatrix}
		\id_{E_0} & 0 \\ 0 & -\id_{E_1}
	\end{pmatrix}.
\end{align*}
On the twist bundle $\overline{\Sigma}^E M \coloneqq \overline{\Sigma} M \otimes E$, we have a $\Z/2\Z$-grading $\iota \otimes \eta$ and an odd $\Cl_1$-action $c_1 \otimes \eta$. The connections on $\overline{\Sigma}M$, $E_0$ and $E_1$ define a connection on $\overline{\Sigma}^E M$, giving rise to a twisted Dirac operator $D^E$ and a twisted $\Cl_1$-linear Dirac-Witten operator $\overline{D}^E = D^E - \frac{1}{2} \tr^g(k) e_0 \cdot \otimes \id_E$, which is again formally self-adjoint, odd and commutes with the $\Cl_1$-action.

The assumption that $E_0$ and $E_1$ agree outside a compact set $K$ leads to even more structure on the twisted bundle. Namely, we choose a smooth cut-off function $\theta \colon M \to [0,1]$ with compact support such that $\theta \equiv 1$ on $K$. Then
\begin{align*}
	T = \iota \otimes \begin{pmatrix}
		0 & -(1 - \theta)i \Psi^{-1} \\ (1 - \theta)i \Psi & 0
	\end{pmatrix}
\end{align*}
is a self-adjoint, odd $\Cl_1$-linear operator that anti-commutes with $\overline{D}^E$ away from $\supp(\upd \theta)$. Hence, for any $\sigma > 0$, the operator $\overline{D}^E + \sigma T$ is again formally self-adjoint, odd and $\Cl_1$-linear and satisfies
\begin{align}
	\left(\overline{D}^E + \sigma T \right)^2 = \left(\overline{D}^E \right)^2 + \sigma^2 T^2 \geq \sigma^2
\end{align}
outside of $\supp(\theta)$. Note that if $M$ itself is compact, then $\theta \equiv 1$ is a possible choice, in which case $T =0$.

We now consider a family of such operators. Here, a technical difficulty arises: The spinor bundle depends on the metric, so the operators act on different bundles. In order to link these bundles we employ a construction similar to the method of generalized cylinders \cite{baer.gauduchon.moroianu:2005} and its futher development, the universal spinor bundle \cite{mueller.nowaczyk:17}. We start by recalling from \cite{baer.gauduchon.moroianu:2005} that a topological spin structure is given by a double covering $P_{\widetilde{\mathrm{GL}_+}(n)} M \to P_{\mathrm{GL}_+(n)} M$ of the principal bundle of positively oriented frames. Moreover, the map associating to a given basis the scalar product, for which this basis is orthonormal, gives rise to an $\SO(n)$-principal bundle $P_{\mathrm{GL}_+(n)} M \to \bigodot_+^2 T^*M$, where $\bigodot_+^2 T^*M$ is the bundle of symmetric positive definite bilinear forms. Now, denoting by $g = (g_t)_{t \in I}$ the family of complete Riemannian metrics obtained by looking at the first component of the path $\gamma \colon I \to \Ini(M)$, we form the pullback squares
\begin{equation*}
	\begin{tikzcd}[column sep=2.5cm]
		{P_{\Spin(n)}(g_\bullet)} \rar \dar & {P_{\widetilde{\mathrm{GL}_+}(n)} M} \dar \\
		{P_{\SO(n)}(g_\bullet)} \rar \dar & {P_{\mathrm{GL}_+(n)} M} \dar \\
		{M \times I} \rar{(x, t) \mapsto g_t(x)} & {\bigodot_+^2 T^*M}.
	\end{tikzcd}
\end{equation*}
Notice that the so-defined principal bundles, the $\SO(n)$-principal bundle $P_{\SO(n)}(g_\bullet) \to M \times I$ and the $\Spin(n)$-principal bundle $P_{\Spin(n)}(g_\bullet) \to M \times I$, are in general just continuous, not smooth, as we assumed the path $\gamma \colon I \to \Ini(M)$ only to be continuous. However, when we restrict to a certain parameter $t \in I$, these give back the (smooth) principal bundles associated to the metric $g_t$, i.\,e.\ the following is a pullback diagram:
\begin{equation*}
	\begin{tikzcd}
		{P_{\Spin(n)}(M, g_t)} \rar \dar & {P_{\Spin(n)}(g_\bullet)} \dar \\
		{P_{\SO(n)}(M, g_t)} \rar \dar & {P_{\SO(n)}(g_\bullet)} \dar \\
		M \rar{x \mapsto (x,t)} & {M \times I}.
	\end{tikzcd}
\end{equation*}

By associating to $P_{\Spin(n)}(g_\bullet)$ the double of the irreducible $\Cl_n$-representation, we obtain a continuous bundle $\overline{\Sigma} g_\bullet \to M \times I$ that restricts for each $t \in I$ to the double spinor bundle $\overline{\Sigma} (M, g_t) \to M$ for the metric $g_t$. Moreover, the twisted bundle $\overline{\Sigma}^E g_\bullet \coloneqq \overline{\Sigma} g_\bullet \otimes p^*E$, where $p \colon M \times I \to M$ is the canonical projection, restricts for fixed $t \in I$ to the twisted bundle $\overline{\Sigma}^E(M, g_t)$ considered above.

A single differential operator on a vector bundle is often best considered as an unbounded operator acting on the $L^2$-space of sections of that bundle.
The corresponding notion for families of operators acting on a family of vector bundles is that of a densely defined operator family on a continuous field of Hilbert spaces. 
This notion was developed by Ebert \cite{ebert:p2018} building on work by Dixmier and Douady \cite{dixmier.douady:63} and in the following we stick to his notation.

We start by constructing a continuous field of Hilbert spaces with $\Cl_1$-structure from the bundle $\overline{\Sigma}^E g_\bullet \to M \times I$. 
Roughly speaking, this consists of the spaces of $L^2$-sections of $\overline{\Sigma}^E g_t \to M$, parametrized over $t \in I$, together with the datum of when a family $(u_t)_{t \in I}$ of $L^2$-sections is continuous.
It will be obtained as completion of an appropriate field of pre-Hilbert spaces:
For each $t \in I$ let
\begin{align*}
	V_t &= C^\infty_c \left(M, \overline{\Sigma}^E(M,g_t) \right)
\intertext{be the space of smooth compactly supported sections and denote by}
	\Lambda &= \left\{ u \in C^0_c \left(M \times I, \overline{\Sigma}^E g_\bullet \right) \,\middle|\, u_{|M \times \{t\}} \in V_t \,\text{for all}\, t \in I \right\} \subseteq \prod_{t \in I} V_t
\end{align*}
the subset of those families of such sections that assemble into a compactly supported continuous section of the bundle $\overline{\Sigma}^E g_\bullet \to M \times I$. 

\begin{Lem}
Together with the $L^2$-scalar product on $V_t$, the pair $((V_t)_{t \in I}, \Lambda)$ defines a continuous field of pre-Hilbert spaces.
\begin{proof}
The first thing to show is that 
\begin{align*}
	I &\lto \R \\	
	t &\lmapsto \int_M \left\langle u(x,t), v(x,t) \right\rangle \upd \mathrm{vol}^{g_t}(x)
\end{align*}
is continuous for $u, v \in \Lambda$. By the definition of $\Lambda$, the functions $(x,t) \mapsto \|u(x,t)\|$ and hence $x \mapsto \max_{t \in I} \|u(x,t)\|$ are continuous and compactly supported, similarly for $v$. Thus, using the Cauchy-Schwarz inequality, the theorem of dominated convergence implies the desired continuity.

Secondly, we have to see that the restriction map $\Lambda \to V_t,\, u \mapsto u_{|M \times \{t\}}$ is dense for all $t \in I$. In fact, we will show surjectivity. So let $u_t \in V_t$ be given. This defines a commutative diagram
\begin{equation*}
	\begin{tikzcd}
		& \overline{\Sigma}^E \dar \\
		M \urar{u_t} \rar{g_t} &  {\bigodot_+^2 T^*M} \arrow[u, bend right, dashed, pos=0.4, "{\tilde{u}}"'],
	\end{tikzcd}
\end{equation*}
whereby $\overline{\Sigma}^E$ denotes the twisted double spinor bundle associated to the spin structure $P_{\widetilde{\mathrm{GL}_+}(n)} M \to \bigodot_+^2 T^*M$. The twist is given by the pull back of $E$ along $\bigodot_+^2 T^*M \to M$. We wish to extend  $u_t$ to a smooth compactly supported section $\tilde{u}$, as this gives rise to a section $u \in \Lambda$ restricting to $u_t$. As $u_t$ is compactly supported, we can turn any smooth extension $\tilde{u}$ into a compactly supported one by multiplying with a suitable cut-off function. We construct $\tilde{u}$ by gluing local pieces using a partition of unity of $\bigodot_+^2 T^*M$: For $x \in M$, let $U \subseteq \bigodot_+^2 T^*M$ be an open neighborhood of $g_t(x)$ with the property that there exists a section $\epsilon$ of $P_{\widetilde{\mathrm{GL}_+}(n)} M_{|U} \to U$. Possibly restricting $U$, we may assume that $\{ g_t(x) \,|\, x \in \pi(U) \} \subseteq U$, where $\pi \colon \bigodot_+^2 T^*M \to M$ is the canonical projection. We now obtain an extension of $u_t$ on $U$ by taking the coefficients with respect to $\epsilon$ in the associated bundle ${\overline{\Sigma}^E}_{|U} \to U$ to be constant along the fibers of $\pi$.
\end{proof}
\end{Lem}

We will denote by $(L^2(\overline{\Sigma}^E g_\bullet), \overline{\Lambda})$ the continuous field of Hilbert spaces obtained as completion of $((V_t)_{t \in I}, \Lambda)$. It is clear that it carries a $\Cl_1$-structure induced by the $\Z/2\Z$-grading $\iota \otimes \eta$ and the Clifford multiplication $c_1 \otimes \eta$. Next, we want to see how $\overline{D}^E + \sigma T$ defines an unbounded operator family on $(L^2(\overline{\Sigma}^E g_\bullet), \overline{\Lambda})$ and establish its main analytic properties with the goal to associate a suitable index to it.

\begin{Lem}
The operators $\overline{D}^E_t + \sigma_t T \colon V_t \to V_t$ for $t \in I$ assemble to a densely defined operator family on $(L^2(\overline{\Sigma}^E g_\bullet), \overline{\Lambda})$ with initial domain $((V_t)_{t \in I}, \Lambda)$.
\begin{proof}
	We have to show that $\overline{D}^E_\bullet + \sigma T$ maps sections $u \in \Lambda$ to sections in $\overline{\Lambda}$. We will show that $(\overline{D}^E_\bullet + \sigma T)u \in \Lambda \subseteq \overline{\Lambda}$. The only thing that is not clear here is that $(\overline{D}^E_\bullet + \sigma T)u$ is a continuous section of $\overline{\Sigma}^E g_\bullet \to M \times I$. This boils down to showing that $D^E_\bullet u$ is a continuous section, where $D^E_\bullet$ is fiberwise the twisted Dirac operator as above. As continuity may be checked locally, we can restrict our attention to an open subset of $M \times I$, where there exists a continuous section $\epsilon$ of $P_{\Spin(n)}(g_\bullet) \to M \times I$. The associated orthonormal frame will be called $(e_1, \ldots, e_n)$. Assuming that $u$ can be written as a tensor product $\Psi \otimes e$ with $\Psi$ a section of $\overline{\Sigma} g_\bullet \to M \times I$ and $e$ a section of $p^* E \to M \times I$ (in general, $u$ will be a sum of such) and expressing $\Psi = [\epsilon, \psi]$, we have
\begin{align*}
	D^E_\bullet u &= \sum_{i =1}^n \left[\epsilon \,,\, e_i \cdot \del_{e_i} \psi + \sum_{j,k =1}^n \frac{1}{2} g_\bullet(\nabla^{g_\bullet}_{e_i} e_j, e_k)\, e_i \cdot e_j \cdot e_k \cdot \psi \right] \otimes e + \sum_{i =1}^n (e_i \cdot \Psi) \otimes \nabla^E_{e_i} e.
\end{align*}
This expression is continuous as we assumed that the family $(g_t)_{t \in I}$ and its derivatives to be uniformly continuous on all compact sets, which particularly implies that the Christoffel symbols are uniformly continuous on compact sets.
\end{proof} 
\end{Lem}

The closure of this operator family will be denoted by $\overline{D}^E_\bullet + \sigma T \colon \mathrm{dom}(\overline{D}^E_\bullet + \sigma T) \to (L^2(\overline{\Sigma}^E g_\bullet), \overline{\Lambda})$.

\begin{Lem} \label{Lem:SelfAdj}
The unbounded operator family $\overline{D}^E_\bullet + \sigma T$ is self-adjoint.
\begin{proof}
By definition, this is the case if and only if the operators $\overline{D}^E_t + \sigma T$ are self-adjoint and regular for each $t \in I$. A sufficient criterion for this is the existence of a coercive, i.~e.\ bounded below and proper, smooth function $h_t \colon M \to \R$, such that the commutator $[\overline{D}^E_t + \sigma_t T,h_t]$ is bounded, cf.\ \cite[Thm. 1.14]{ebert:p2018}.

For some fixed $x_0 \in M$, we consider the distance function $d_t(x) = \mathrm{dist}^{g_t}(x_0,x)$. This is bounded below by $0$. As $g_t$ is complete, by the theorem of Hopf-Rinow, $d_t$ is proper. Yet $d_t$ is not smooth, only a Lipschitz function, with Lipschitz constant $1$. The desired function $h_t$ can now be obtained by suitably smoothing $d_t$ out. For example, using the main theorem of \cite{azagra.ferrera.lopez-mesas.rangel:2007}, we obtain the existence of a smooth function $h_t$ with $\sup_{x \in M} |h_t(x) - d_t(x)| \leq \epsilon$ and Lipschitz constant $1 + \epsilon$ (for any $\epsilon > 0$). The first condition ensures that $h_t$ is also bounded below and proper, whereas the second one implies that $[\overline{D}^E_t + \sigma_t T,h_t]$ is bounded by $1 + \epsilon$, as the principal symbol of $\overline{D}^E_t + \sigma_t T$ is given by Clifford multiplication. 
\end{proof}
\end{Lem}

\begin{Prop} \label{Prop:FredProp}
The self-adjoint unbounded operator family $\overline{D}^E_\bullet + \sigma T$ is a Fredholm family.
\begin{proof}
The first thing to note is that $\overline{D}^E_\bullet + \sigma T$ arises as closure of a formally self-adjoint \emph{elliptic} differential operator of order $1$ on the bundle $\overline{\Sigma}^E g_\bullet \to M \times I$. In view of \cite[Thm. 2.41]{ebert:p2018}, the statement is basically a consequence of the fact that $(\overline{D}^E_\bullet + \sigma T)^2 \geq \sigma^2$ outside the compact set $\supp(\theta) \times I$. However, we are not precisely in the setting of Ebert's article. Namely, the bundle $\overline{\Sigma}^E g_\bullet \to M \times I$ is only continuous and not smooth; but this lower regularity does not affect the proofs. Moreover, we have not yet established the existence of a smooth coercive function $h \colon M \times I \to \R$ such that $[\overline{D}^E_\bullet + \sigma T,h]$ is bounded -- and we will not do so.

Instead, we note that $h$ serves only two purposes. Firstly, it (again) shows that the operator family is self-adjoint, as the functions $h(-,t)$ can play the role of the $h_t$ above. Secondly, it serves as a basis for constructing a compactly supported smooth function $f \colon M \times I \to \R$ such that $(\overline{D}^E_\bullet + \sigma T)^2 + f^2 \geq \sigma^2$ everywhere on $M \times I$ and $\|[\overline{D}^E_\bullet + \sigma T,f]\| \leq \frac{\sigma^2}{2}$, which is needed in the proof of Fredholmness. So we may just construct such a function $f$ directly.

For any $t \in I$, let $R_t$ be chosen such that $\supp(\theta) \subseteq B_{R_t}^{g_t}(x_0)$ and $h_t$ a function as above (for $\epsilon \leq \frac{1}{3}$). Furthermore, we choose a smooth cut-off function $\Psi \colon \R \to [0,1]$ with $\Psi(r) = \sigma$ for $r \leq 1$ and $|\Psi'(r)| \leq \frac{\sigma^2}{3}$ for all $r \in \R$. Denote by $L$ a number such that $\Psi \equiv 0$ on $[L, \infty)$. Now, let $f_t(x) = \Psi(h_t(x) - R_t)$. Note that $f_t \equiv \sigma$ on $\supp(\theta)$, as $h_t(x) \leq d_t(x)+ \epsilon < R_t + 1$ for all $x \in \supp(\theta)$.

Continuity of $(g_s)_{s \in I}$ allows us to choose $\delta_t > 0$ such that $\|g_s^{-1}\|_{g_t} \leq \left(\frac{9}{8}\right)^2$ for all $s \in U_t \coloneqq (t-\delta_t, t+\delta_t)$ on $\overline{B_{R_t+L}^{g_t}(x_0)}$. Using
$\| \upd f_t\|_{g_s}^2 \leq \|g_s^{-1}\|_{g_t} \|\upd f_t\|_{g_t}^2$, we obtain
\begin{align} \label{eq:LemFredIneq} \nonumber
	\|[\overline{D}_s^E + \sigma T,f_t]\| &= \|\upd f_t\|_{g_s} \\
		&\leq \sqrt{\|g_s^{-1}\|_{g_t}}  \|\Psi'\|_\infty \|\upd h_t\|_{g_t} \\ \nonumber
		&\leq \frac{9}{8} \cdot \frac{\sigma^2}{3} \cdot (1+\epsilon)
		\leq \frac{\sigma^2}{2}
\end{align}
for all $s \in U_t$.

Now, there exists a finite collection $t_1, \ldots, t_n$ such that $U_{t_1}, \ldots, U_{t_n}$ cover $I$ and a smooth partition of unity $\psi_1, \ldots \psi_n$ subordinate to this open cover. We define $f(x,t)= \sum_{i=1}^n \psi_i(t) f_t(x)$. Then $f \equiv \sigma$ on $\supp(\theta)$, which implies $(\overline{D}^E_\bullet + \sigma T)^2 + f^2 \geq \sigma^2$ everywhere on $M \times I$. Moreover, the second property of $f$ immediately follows from \eqref{eq:LemFredIneq}.
\end{proof}
\end{Prop}

It is important to know when the Fredholm family $\overline{D}^E_\bullet + \sigma T$ is invertible. The next lemma provides a criterion for this.
\begin{Lem} \label{Lem:Invertibility}
Let $A \subseteq I$ with
\begin{align} \label{eq:compcond}
	\inf_{t \in A,\, x \in M} \left(\rho_t - \|j_t\|_{g_t} - 2\|\mathcal{R}^E_t\|\right) &> 0.
\end{align}
Then there is a $\sigma' > 0$ such that $\overline{D}^E_\bullet + \sigma T$ is invertible over $A$ for all $0 < \sigma < \sigma'$. Here, $\mathcal{R}^E_t$ denotes the curvature endomorphism defined by $\mathcal{R}^E_t(\phi \otimes e)= \sum_{i<j} e_i \cdot  e_j \cdot \phi \otimes R^E(e_i,e_j)e$ for $\phi \otimes e \in (\Sigma_p M \oplus \Sigma_p M) \otimes_\C E_p$ and an orthonormal basis $(e_1, \ldots, e_n)$ of $T_pM$, $p \in M$, with respect to $g_t$. 
\begin{proof}
	We first consider the situation for some fixed metric $g$. The twisted Dirac-Witten operator associated to $g$ satisfies the following Schrödinger-Lichnerowicz type formula, the proof of which is deferred to the appendix:
	\begin{align*}
	\left(\overline{D}^E\right)^2 \psi &= (\overline{\nabla}^E)^* \overline{\nabla}^E \psi + \frac{1}{2} (\rho - e_0 \cdot j^\sharp \cdot) \psi + \mathcal{R}^E \psi.
	\end{align*}
	Here, $\overline{\nabla}^E$ denotes the connection on $\overline{\Sigma}^E M$ induced by $\overline{\nabla}$ and the connection on $E$, and the star indicates the formal adjoint. Together with
	\begin{align*}
		\left(\overline{D}^E + \sigma T\right)^2 \psi &= \left(\overline{D}^E \right)^2 \psi + \sigma \left(\overline{D}^E T + T \overline{D}^E \right) \psi + \sigma^2 (1-\theta)^2 \psi
	\end{align*}
	this implies 
	\begin{align*}
		\left(\left(\overline{D}^E + \sigma T\right)^2 \psi, \psi \right)_{L^2}
		&\geq \left(\left(\overline{D}^E\right)^2 \psi, \psi \right)_{L^2} -\sigma \|\upd \theta\|_g \|\psi\|_{L^2}^2  + \sigma^2 (1-\theta)^2 \|\psi\|_{L^2}^2   \\
		&\geq \|\overline{\nabla}^E \psi \|_{L^2}^2 + \left(\frac{1}{2}(\rho - \|j\|_g)-\|\mathcal{R}^E\| \right) \|\psi\|_{L^2}^2 -\sigma \|\upd \theta\|_g \|\psi\|_{L^2}^2 \\
		 &\geq \left(\frac{1}{2}(\rho - \|j\|_g)-\|\mathcal{R}^E\|\right) \|\psi\|_{L^2}^2 -\sigma \|\upd \theta\|_g \|\psi\|_{L^2}^2.
	\end{align*}
	for any compactly supported smooth section $\psi$ (norms without subscript $L^2$ denote pointwise norms).
	
	As $\|d\theta\|_{g_t}$ is a continuous compactly supported function on $M \times I$ and \eqref{eq:compcond}	holds, we may choose $\sigma'> 0$ such that
	\begin{align*}
		\inf_{t \in A,\, x \in M} \left(\rho_t - \|j_t\|_{g_t} - 2\|\mathcal{R}^E_t\|\right) \geq \sigma' \sup_{t \in A,\, x \in M} \|d\theta\|_{g_t}.
	\end{align*}
	Then for all $0 < \sigma < \sigma'$, there exists some constant $c>0$ with $\left(\overline{D}^E_t + \sigma T\right)^2 \geq c$	for all $t \in A$.
	From this, the statement follows immediately (cf.\ \cite{ebert:p2018}[Prop.\ 1.21, Lem.\ 2.6]).
\end{proof}
\end{Lem}

\begin{Prop} \label{Prop:DefInd}
The operator family $\overline{D}^E_\bullet + \sigma T$ is odd with respect to $\iota \otimes \eta$ and $\Cl_1$-linear with respect to $c_1 \otimes \eta$. For suitably small $\sigma$, it defines an element
\begin{align*}
	\left[(L^2(\overline{\Sigma}^E g_\bullet), \overline{\Lambda}),\; \iota \otimes \eta,\; -i\iota c_1 \otimes \id_E,\; \overline{D}^E_\bullet + \sigma T \right] &\in \mathrm{K}^1(I, \del I),
\end{align*}
that is independent of the choices of $K, \theta, \Psi$ and $\sigma$ (as long as they fulfill the assumed requirements) and depends only on the relative homotopy class of $\gamma \colon (I, \del I) \to (\Ini(M), \Ini^{> \mathcal{R}^E +c}(M))$. Moreover, this class is zero if $\gamma$ is homotopic to a path $I \to \Ini^{> \mathcal{R}^E +c}(M)$.
\end{Prop}
Note that, as $\Ini(M)$ is convex, the homotopy class of $\gamma$ just depends on its endpoints. 

\begin{Def}
For $(g_{-1},k_{-1}), (g_{1},k_{1}) \in \Ini^{> \mathcal{R}^E +c}(M)$, their \emph{$E$-relative index difference} $\idiff^E((g_{-1},k_{-1}), (g_{1},k_{1})) \in \mathrm{K}^1(I, \del I)$ is the class defined in \cref{Prop:DefInd} using some path $\gamma \colon (I, \del I) \to (\Ini(M), \Ini^{> \mathcal{R}^E +c}(M))$ connecting these two pairs.
\end{Def}

If $M$ is compact, then $E_0$ and $E_1$ need not be isomorphic anywhere. In particular, we may take $E_0 = \C$ and $E_1 = 0$, which gives the untwisted index difference.

\begin{Def}
If $M$ is compact, the \emph{index difference} of $(g_{-1},k_{-1})$ and $(g_{1},k_{1}) \in \Ini^{+}(M)$  is $\idiff((g_{-1},k_{-1}), (g_{1},k_{1})) \coloneqq \idiff^E((g_{-1},k_{-1}), (g_{1},k_{1}))\in \mathrm{K}^1(I, \del I)$ for the trivial bundles $E_0 = \underline{\C}$ and $E_1 = \underline{0}$.
\end{Def} 

\begin{proof}
That $\overline{D}^E_\bullet + \sigma T$ is odd and $\Cl_1$-linear follows from the fact that this holds for $\overline{D}^E_t + \sigma T$ for all $t \in I$. As $\overline{D}^E_\bullet + \sigma T$ is, moreover, an unbounded Fredholm family that is by \cref{Lem:Invertibility} invertible over $\del I$, we get an element in the K-theory of $(I,\del I)$. Note that the K-theoretic model described in \cite{ebert:p2018}[Ch.\ 3], which we are using here, requires a $\Cl_1$-antilinear operator rather than a $\Cl_1$-linear one. But this is no problem as a $\Cl_1$-linear operator is $\Cl_1$-antilinear with respect to the $\Cl_1$-structure defined by $i(\iota \otimes \eta)(c_1 \otimes \eta) = -i\iota c_1 \otimes \id$.

We now show independence of the choices starting with $\sigma$. Let $\sigma_0 > 0$ and $\sigma_1 >0 $ be two admissible values, i.~e.\ smaller than $\sigma'$ from \cref{Lem:Invertibility} applied to $A = \del I$. We consider the pullback of the $\left((L^2(\overline{\Sigma}^E g_\bullet), \overline{\Lambda}),\; \iota \otimes \eta,\; -i\iota c_1 \otimes \id \right)$ along the canonical projection $I \times [0,1] \to I$. This continuous field of $\Z/2\Z$-graded Hilbert spaces with $\Cl_1$-structure carries the odd $\Cl_1$-antilinear Fredholm family $\overline{D}^E_t + ((1-s)\sigma_0 + s\sigma_1)T$ for $t \in I$ and $s \in [0,1]$, which is invertible over $\del I \times [0,1]$. Thus, we have a concordance between the cycles $\left((L^2(\overline{\Sigma}^E g_\bullet), \overline{\Lambda}),\; \iota \otimes \eta,\; -i\iota c_1 \otimes \id,\; \overline{D}^E_\bullet + \sigma T \right)$ for $\sigma = \sigma_0$ and $\sigma = \sigma_1$.

Independence of $\theta$, $K$ and $\Psi$ are slightly connected, as we have to have $\theta \equiv 1$ on $K$ and $\Psi$ is defined on the complement of $K$. Given two such triples $(\theta_0, K_0, \Psi_0)$ and $(\theta_1, K_1, \Psi_1)$, we first show that we can first replace the $\theta_0$ by some $\theta$ with $\theta \equiv 1$ on $K = K_0 \cup K_1$ without changing the $\mathrm{K}^1$-class. Then noting that $(\theta, K_0, \Psi_0)$ and $(\theta, K, {\Psi_0}_{|M \setminus K})$ even define the same operator family, it just remains to show that the $\mathrm{K}^1$-class is independent of $\Psi$ for fixed $\theta$ and $K$.

Concering the replacement of $\theta_0$ by $\theta$ (similarly for $\theta_1$ by $\theta$), we use the same argumentation as for $\sigma$ with the difference that this time the operator family is given by $\overline{D}^E_t + \sigma T_s$ with
\begin{align*}
	T_s =  \iota \otimes \begin{pmatrix}
		0 & -(1 - (1-s)\theta_0 - s\theta)i \Psi_0^{-1} \\ (1 - (1-s)\theta_0 -s \theta)i \Psi_0 & 0
	\end{pmatrix}.
\end{align*}

For changing ${\Psi_0}_{|M \setminus K}$ to ${\Psi_1}_{|M \setminus K}$, we note that the requirement that these bundle isomorphisms preserve hermitian metric and connection implies that they differ by a single element of $U(k)$ on every connected component of $M \setminus K$, where $k$ is the rank of $E_0$ (and $E_1$). As $U(k)$ is connected, there exists a homotopy $(\Psi_s)_{s \in [0,1]}$ connecting these. Again, the operator family $\overline{D}^E_t + \sigma T_s$ defines a concordance, where this time
\begin{align*}
	T_s =  \iota \otimes \begin{pmatrix}
		0 & -(1 - \theta)i \Psi_s^{-1} \\ (1 - \theta)i \Psi_s & 0
	\end{pmatrix}.
\end{align*}

Now assume that we are given a homotopy $H \colon (I \times [0,1], \del I \times [0,1]) \to (\Ini(M), \Ini^{> \mathcal{R}^E +c}(M))$ between $\gamma_0 = H(-,0)$ and $\gamma_1=H(-,1)$. In this case we can construct a concordance between the $\mathrm{K}^1$-cycles associated to $\gamma_0$ and $\gamma_1$ as follows. Let $(g_{t,s})_{t \in I,\, s \in [0,1]}$ be the first components of the pairs $(H(t,s))_{t \in I,\, s \in [0,1]}$. Similarly to before, we may form the pullback
\begin{equation*}
	\begin{tikzcd}[column sep=2.5cm]
		{P_{\Spin(n)}(g_{\bullet,\bullet})} \rar \dar & {P_{\widetilde{\mathrm{GL}_+}(n)} M} \dar \\
		{P_{\SO(n)}(g_{\bullet,\bullet})} \rar \dar & {P_{\mathrm{GL}_+(n)} M} \dar \\
		{M \times I \times [0,1]} \rar{(x, t, s) \mapsto g_{t,s}(x)} & {\bigodot_+^2 T^*M},
	\end{tikzcd}
\end{equation*}
obtain a bundle $\overline{\Sigma}^E g_{\bullet,\bullet} \to M \times I \times [0,1]$ and a continuous field of $\Z/2\Z$-graded Hilbert spaces $(L^2(\overline{\Sigma}^E g_{\bullet,\bullet}),\overline{\Lambda})$ with $\Cl_1$-structure. The operators $(\overline{D}^E_{t,s} + \sigma T)_{t \in I,\, s \in [0,1]}$ constitute an unbounded Fredholm family on this continuous field of Hilbert spaces. For suitably small $\sigma$, by an analogous statement to \cref{Lem:Invertibility}, this is invertible over $\del I \times [0,1]$. Together, these provide the required concordance.

For the last statement assume that the image of $\gamma$ is contained in $\Ini^{> \mathcal{R}^E +c}$. In this case, $\overline{D}^E_\bullet + \sigma T$ is invertible on all of $I$ by \cref{Lem:Invertibility}, as long as $\sigma$ is chosen suitably small. Therefore, the $\mathrm{K}^1$-class under consideration is in the image of the restriction homomorphism $\mathrm{K}^1(I,I) \to \mathrm{K}^1(I, \del I)$. But as $\mathrm{K}^1(I,I)=0$, this class has to be zero.
\end{proof}

\section{An index theorem for the twisted index difference}
The purpose of this section is to calculate the $E$-relative index difference for between pairs of the form $(g, \tau g)$ and $(g, -\tau g)$ for some $\tau > 0$. This is done in two steps. The first one is to express the relative index difference as relative index of a suitable Dirac type operator. In the second step, a relative index theorem identifies this index with a twisted version of the $\hat{A}$-genus.

The first thing we realize is the following: If for $(g_{-1},k_{-1})$ and $(g_1,k_1)$ there exists some $\sigma' > 0$ such that $\overline{D}^E_\bullet + \sigma T$ is invertible for all $0 < \sigma < \sigma'$, then it makes sense to speak of their $E$-relative index difference, even if they are not contained in some $\Ini^{>\mathcal{R}^E +c}(M)$.

Let $g$ be some complete metric on $M$. Such a metric always exists by a classical result of Greene \cite{greene:1978}, it is even the case that every conformal class contains such a metric by \cite{mueller.nardmann:2015}. For $\tau > 0$, we consider the pairs $(g, \tau g)$ and $(g, -\tau g)$. Using
\begin{align*}
	\left(\overline{D}^E \right)^2 
	&= \left(D^E \pm \frac{1}{2} \tau e_0 \cdot \otimes \id_E \right)^2 
	= \left(D^E \right)^2 + \frac{n^2}{4} \tau^2
\intertext{we obtain}
		\left(\overline{D}^E + \sigma T\right)^2
		&= \left(\overline{D}^E \right)^2 + \sigma \left(\overline{D}^E T + T \overline{D}^E \right) + \sigma^2 (1-\theta)^2 \\
		&\geq \left(D^E \right)^2 + \frac{n^2}{4} \tau^2 - \sigma \|\upd \theta\|_g,
\end{align*}
which shows that for these pairs $\overline{D}^E + \sigma T$ is invertible for sufficiently small $\sigma > 0$. Thus it makes sense to speak of $\idiff^E((g,-\tau g), (g, \tau g)) \in \mathrm{K}^1(I, \del I)$.

We denote by $D^E_0$ the Dirac operator on $\Sigma^E M \coloneqq \Sigma M \otimes E$ for the metric $g$. With
\begin{align*}
	T_0 = \omega \otimes \begin{pmatrix}
		0 & -(1 - \theta)i \Psi^{-1} \\ (1 - \theta)i \Psi & 0
	\end{pmatrix},
\end{align*}
this gives rise to a Fredholm operator $D^E_0 + \sigma T_0$, for $\sigma > 0$, which is odd with respect to the $\Z/2\Z$-grading $\omega \otimes \eta$. The proof of Fredholmness uses that $g$ is complete. It is similar to \cref{Lem:SelfAdj,Prop:FredProp}, but much less delicate.

\begin{Prop}
The isomorphism $\mathrm{K}^1(I, \del I) \cong \mathrm{K}^0(\pt) \cong \Z$, induced by the Bott map, maps the class $\idiff^E((g,-\tau g), (g, \tau g))$ to $\ind(D^E_0 + \sigma T_0)$.
\begin{proof}
	We start from the class $[L^2(\Sigma^E M),\; \omega \otimes \eta,\; D^E_0 + \sigma T_0] \in \mathrm{K}^0(\pt)$ corresponding to the integer value $\ind(D^E_0 + \sigma T_0)$. Using the conventions of Ebert \cite{ebert:p2018}, the Bott map $\mathrm{K}^0(\pt) \overset{\sim}{\to} \mathrm{K}^1(\R, \R\setminus \{0\})$ sends this class to
\begin{align*}
	\left[\C^2 \otimes p^*L^2(\Sigma^E M),
		\begin{pmatrix}
			\omega \otimes \eta & 0 \\ 0 & -\omega \otimes \eta
		\end{pmatrix},
		\begin{pmatrix}
			0 & -\omega \otimes \eta \\ \omega \otimes \eta & 0
		\end{pmatrix},
		 \begin{pmatrix}
			D^E_0 + \sigma T_0 & t \omega \otimes \eta \\ 
			t \omega \otimes \eta & D^E_0 + \sigma T_0
		\end{pmatrix}
	\right],
\end{align*}
where $p \colon \R \to \pt$ is the projection and $t$ is the $\R$-coordinate. As the inclusion $(I, \del I) \to (\R, \R \setminus \{0\})$ induces an isomorphism in $K$-theory, the same formula defines the corresponding element in $\mathrm{K}^1(I,\del I)$; now assuming $p \colon I \to \pt$ and $t \in I$.

Note that we may identify $\left(L^2(\overline{\Sigma}^E g_\bullet\right), \overline{\Lambda}) = \C^2 \otimes p^* L^2(\Sigma^E M)$ for the constant family $g_t = g$. Using the automorphism of $\left(L^2(\overline{\Sigma}^E g_\bullet), \overline{\Lambda}\right)$ given by
\begin{align*}
	\begin{pmatrix}
		\id_{\Sigma M} & 0 \\ 0 & \id_{\Sigma M} 
	\end{pmatrix}
	\otimes
	\begin{pmatrix}
		\id_{E_0} & 0 \\ 0 & 0
	\end{pmatrix} +
		\begin{pmatrix}
		\id_{\Sigma M} & 0 \\ 0 & -\id_{\Sigma M} 
	\end{pmatrix}
	\otimes
	\begin{pmatrix}
		0 & 0 \\ 0 & \id_{E_1}
	\end{pmatrix},
\end{align*}
the $\mathrm{K}^1$-class translates into
\begin{align*}
	\left[\left(L^2(\overline{\Sigma}^E g_\bullet), \overline{\Lambda}\right),
		\begin{pmatrix}
			\omega & 0 \\ 0 & -\omega
		\end{pmatrix} \otimes \eta ,
		\begin{pmatrix}
			0 & -\omega \\ \omega  & 0
		\end{pmatrix} \otimes \id_E ,
		 \begin{pmatrix}
			D^E_0 + \sigma T_0 & t \omega \otimes \id_E \\ 
			t \omega \otimes \id_E & D^E_0 - \sigma T_0
		\end{pmatrix}
	\right].
\end{align*}
Applying furthermore the automorphism
\begin{align*}
	\frac{1}{\sqrt 2} \left(\begin{pmatrix}
	\id_{\Sigma M} & 0 \\ 0 & \id_{\Sigma M}
\end{pmatrix}	 + \begin{pmatrix}
		0 & \id_{\Sigma M} \\ -\id_{\Sigma M} & 0
\end{pmatrix} \right) \otimes \id_E,
\end{align*}
this gets
\begin{gather*}
	\left[\left(L^2(\overline{\Sigma}^E g_\bullet), \overline{\Lambda}\right),
		\begin{pmatrix}
			0 & \omega \\ \omega & 0
		\end{pmatrix} \otimes \eta ,
		\begin{pmatrix}
			0 & -\omega \\ \omega  & 0
		\end{pmatrix} \otimes \id_E,
		 \begin{pmatrix}
			D^E_0 - t \omega \otimes \id_E & \sigma T_0 \\ 
			\sigma T_0 & D^E_0 + t \omega \otimes \id_E 
		\end{pmatrix}
	\right] \\
	= 
	\left[\left(L^2(\overline{\Sigma}^E g_\bullet), \overline{\Lambda}\right),\;
		\iota \otimes \eta ,\;
		-i\iota c_1 \otimes \id_E,\;
		D^E - t e_0 \cdot \otimes \id_E + \sigma T
	\right].
\end{gather*}
Bearing in mind that $\overline{D}^E = D^E - \frac{1}{2} \tr(k) e_0 \cdot \otimes \id_E$, this is the class defined by the Dirac-Witten operators associated to the straight unit speed path from $(g,-2g)$ to $(g, 2g)$. For general $\tau > 0$, the result follows by rescaling, e.~g.\  replacing the inclusion $(I, \del I) \to (\R, \R \setminus \{0\})$ above with the map $t \mapsto \frac{\tau}{2} t$.
\end{proof}
\end{Prop}

It now remains to determine $\ind(D^E_0 + \sigma T_0)$. This is done by the relative index theorem going back to \cite{gromov.lawson:1983}.
\begin{Satz}[Relative Index Theorem]
\begin{align*}
	\ind(D^E_0 + \sigma T_0) = \int_M \hat{A}(TM) \wedge (\mathrm{ch}(E_1)-\mathrm{ch}(E_0)) \eqqcolon \hat{A}(M,E).
\end{align*}
\begin{proof}
	Although probably well-known, there seems not to be an easily citable reference matching the setup here.
	We therefore provide a proof using cut-and-paste-techniques from \cite{baer.ballmann:2011}.
	
	First, we note that the $\Z/2\Z$-graded index of $D^E_0 + \sigma T_0$ is by definition just the usual Fredholm index of its "positive" part $(D^E_0)^+ + \sigma T_0^+ \colon \Gamma((\overline{\Sigma}^E M)^+) \to \Gamma((\overline{\Sigma}^E M)^-)$ mapping from positive to negative half-spinors. The operator $(D^E_0)^+ + \sigma T_0^+$ is of Dirac type, so we may use the decomposition theorem from \cite{baer.ballmann:2011}.
	
	In order to do so, let $M = M_1 \cup M_2$ be a decomposition into two smooth manifolds with boundary $\del M_1 = \del M_2$ such that $M_1$ is compact and $\supp \theta \subseteq M_1$. Let $B \subseteq H^{\frac12}(\del M_1, (\overline{\Sigma}^E M)^+)$ be an elliptic boundary condition and denote by $B^\perp$ its $L^2$-orthogonal complement. For example, $B$ could be Atiyah-Patodi-Singer boundary conditions. The decomposition theorem states that
	\begin{align*}
		\Ind\left((D^E_0)^+ + \sigma T_0^+\right) &= \Ind\left(((D^E_0)^+ + \sigma T_0^+)_{|M_1,\,B}\right) + \Ind\left(((D^E_0)^+ + \sigma T_0^+)_{|M_2,\,B^\perp}\right).
	\end{align*}
	
	The first summand computes to
	\begin{equation} \label{eq:IndInt}
	\begin{aligned} 
		\Ind\left(((D^E_0)^+ + \sigma T_0^+)_{|M_1,\,B}\right)
			&= \Ind\left((D^E_0)^+_{|M_1,\,B}\right) \\
			&= \Ind\left((D^{E_0}_0)^+_{|M_1,\,B}\right) - \Ind\left((D^{E_1}_0)^+_{|M_1,\,B}\right) \\
			&= \int_M \hat{A} \wedge (\mathrm{ch}(E_1)-\mathrm{ch}(E_0)).
	\end{aligned}	
	\end{equation}
	Here, in the first step, we used the homotopy $[0,1] \ni t \mapsto (D^E_0)^+ + t \sigma T_0^+$. The second step follows from the decomposition
	\begin{align*}
		(D^E_0)^+ &= \begin{pmatrix}
			(D^{E_0}_0)^+ & 0 \\ 0 & (D^{E_1}_0)^-
		\end{pmatrix}
		=\begin{pmatrix}
			(D^{E_0}_0)^+ & 0 \\ 0 & \left((D^{E_1}_0)^+ \right)^*
		\end{pmatrix}.
	\end{align*}
	The last step is (contained in the proof of) the relative index theorem \cite[Thm. 1.21]{baer.ballmann:2011}.
	
	It remains to show that the second summand is zero. To see this, we note that there exists a bundle $\widetilde{E}_1 \to M$ admitting a metric and connection preserving isomorphism $\widetilde{\Psi} \colon E_0 \to \widetilde{E}_1$, such that $\widetilde{E}_{1|M_2} = {E_1}_{|M_2}$ and $\widetilde{\Psi}_{|M_2} = \Psi_{|M_2}$. For instance, such a bundle can be obtained by gluing $E_{0|M_1}$ and ${E_1}_{|M_2}$. Similarly as before, we denote by $D_0^{\widetilde{E}}$ the Dirac operator on $\Sigma^{\widetilde{E}} M$ for $\widetilde{E} = E_0 \oplus \widetilde{E}_1$ and define
	\begin{align*}
	\widetilde{T}_0 = \omega \otimes \begin{pmatrix}
		0 & -i \widetilde{\Psi}^{-1} \\ i \widetilde{\Psi} & 0
	\end{pmatrix}.
	\end{align*}
	Notice, that we are allowed to take $\theta \equiv 0$, as $\widetilde{\Psi}$ is defined on all of $M$. Again, we have a decomposition
	\begin{align*}
		\Ind\left((D^{\widetilde{E}}_0)^+ + \sigma \widetilde{T}_0^+\right) &= \Ind\left(((D^{\widetilde{E}}_0)^+ + \sigma \widetilde{T}_0^+)_{|M_1,\,B}\right) + \Ind\left(((D^{\widetilde{E}}_0)^+ + \sigma \widetilde{T}_0^+)_{|M_2,\,B^\perp}\right).
	\end{align*}
	In this case, the calculation \eqref{eq:IndInt} shows that the first summand is zero, as $\mathrm{ch}(\widetilde{E}_1) = \mathrm{ch}(E_0)$. The second summand is the second summand from above as the bundles and operators are equal on $M_2$. Thus, we obtain
	\begin{align*}
		\Ind\left(((D^E_0)^+ + \sigma T_0^+)_{|M_2,\,B^\perp}\right) 
		&= \Ind\left(((D^{\widetilde{E}}_0)^+ + \sigma \widetilde{T}_0^+)_{|M_2,\,B^\perp}\right) \\
		&= \Ind\left((D^{\widetilde{E}}_0)^+ + \sigma \widetilde{T}_0^+\right) = 0,
	\end{align*}
	where we used that $D^{\widetilde{E}}_0 + \sigma \widetilde{T}_0$ is invertible as $(D^{\widetilde{E}}_0 + \sigma \widetilde{T}_0)^2 \geq \sigma^2 > 0$.	
\end{proof}
\end{Satz}

\begin{Kor}[Relative index theorem for the twisted index difference] \label{Cor:RelIndThm}
The isomorphism $\mathrm{K}^1(I, \del I) \cong \mathrm{K}^0(\pt) \cong \Z$ sends $\idiff^E((g,-\tau g), (g, \tau g))$ to $\hat{A}(M,E)$. If $M$ is compact, $\idiff((g,-\tau g), (g, \tau g))$ is sent to $\hat{A}(M)$.
\end{Kor}

The relative index theorem allows to obtain an obstruction to the path-connectedness of the space of initial data sets that strictly satisfy the dominant energy condition. The following corollary illustrates the general strategy and turns out to be a special case of the enlargeability obstruction \cref{Thm:EnlObst} that we discuss in the remaining section.

\begin{Kor} \label{Cor:Ahat}
Let $M$ be a compact manifold, $g$ a metric on $M$ and $\tau > 0$ be chosen such that $(g,-\tau g), (g, \tau g) \in \DEC(M)$.
	If $M$ is spin and $\hat{A}(M) \neq 0$, then $(g,-\tau g)$ and  $(g, \tau g)$ belong to different path-components of $\DEC(M)$.
	\begin{proof}
		If there were a path $\gamma \colon I \to \DEC(M)$ from $(g,-\tau g)$ to $(g, \tau g)$, then, by \cref{Prop:DefInd}, $\idiff((g,-\tau g), (g, \tau g))$ would be zero. But by \cref{Cor:RelIndThm} it is mapped to $\hat{A}(M) \neq 0$.
	\end{proof}
\end{Kor}

\begin{Bem}
	The statement of \cref{Cor:Ahat} also follows from the main result in \cite{gloeckle:p2019}. This is due to the fact that, in even dimension $n$, Hitchin's $\alpha$-index is mapped to the $\hat A$-genus under the complexification map $\KO^{-n}(\pt) \to \mathrm{K}^{-n}(\pt) \cong \Z$. In fact, also the proof is the same as complexification turns the $\alpha$-difference defined in \cite{gloeckle:p2019} into the (untwisted) index-difference considered in this article, up to Bott periodicity.
\end{Bem}

\section{Enlargeability obstruction for initial data sets}
Gromov-Lawson's enlargeability obstruction gives a major source of examples of manifolds that do not admit a positive scalar curvature metric. In this section, we prove that enlargeability is also an obstruction to path-connectedness of the space of initial data sets strictly satisfying the dominant energy condition. There are many versions of enlargeability.
In what follows, we will always understand enlargeability in the sense of $\hat{A}$-area-enlargeability:

\begin{Def} \label{Def:Enl}
	A smooth map $f \colon (M,g) \to (N,h)$ between Riemannian manifolds is \emph{$\epsilon$-area-contracting} for some $\epsilon > 0$ if the induced map $f_* \colon \Lambda^2 TM \to \Lambda^2 TN$ satisfies $\|f_*\| \leq \epsilon$.
	A compact Riemannian manifold $(M,g)$ of dimension $n$ is called \emph{area-enlargeable in dimension $k$} if for all $\epsilon > 0$ there exists a Riemannian covering $(M^\prime,g^\prime) \to (M,g)$ admitting an $\epsilon$-area-contracting map $(M^\prime,g^\prime) \to (S^k,g_{Std})$ that is constant outside a compact set and of non-zero $\hat A$-degree. It is called \emph{$\hat A$-area-enlargeable} if it is area-enlargeable in some dimension $k$.
\end{Def}

Recall that the $\hat A$-degree of a smooth map $f \colon X \to Y$, where $Y$ is compact and connected and $f$ is constant outside a compact set, may be defined by the requirement that $\int_X \hat A(TX) \wedge f^*(\omega) = \A(f) \int_Y \omega$ for all top dimensional forms $\omega \in \Omega^{\dim(Y)}(Y)$. If $Y$ is non-connected, there is one such number for every connected component of $Y$ and the $\hat A$-degree is the vector consisting of these. From this definition we see that it can only be non-zero if $\dim(Y) \leq \dim(X)$ and $\dim(Y) \equiv \dim(X) \mod 4$. The $\hat A$-degree can be thought of as interpolating between the following two special cases: If $\dim(Y) = \dim(X)$, the $\hat A$-degree is just the usual degree. If $\dim(Y)=0$ and $Y$ is connected, it is the $\hat A$-genus $\hat A(X)$ of $X$. 

Although the definition of enlargeability uses a Riemannian metric, the property itself is independent of this choice. Manifolds that are enlargeable in dimension $0$ are precisely the ones having non-zero $\hat A$-genus. Another main example is the torus $T^n = \R^n/\Z^n$, which is enlargeable in the top dimesion $n$. Furthermore, every compact manifold that admits a metric of non-positive sectional curvature is enlargeable (in the top dimension) by the Cartan-Hadarmard theorem. If $M$ is enlargeable then the direct sum $M \# N$ with another manifold $N$ is again enlargeable. Furthermore, for an enlargeable manifold $M$ the product $M \times S^1$ with a circle is again enlargeable. This, and much more, is discussed in great detail in \cite[Sec.~IV.5]{lawson.michelsohn:89}.

The proof that enlargeable spin manifolds do not admit psc metrics \cite[Thm.~5.21]{gromov.lawson:1983} can be split into two parts. The first part consists of using the enlargeability condition to construct a suitable family of complex vector bundles over coverings of the manifold:
\begin{Satz}[Gromov-Lawson] \label{Thm:GL}
	Let $(M,g)$ be an $\hat A$-area-enlargeable manifold of even dimension. Then there exists a sequence of coverings $M_i \to M$ and $\Z/2\Z$-graded hermitian vector bundles $E_i = E_i^{(0)} \oplus E_i^{(1)} \to M_i$ with compatible connection, such that
	\begin{itemize}
		\item for all $i \in \N$ the bundles $E_i^{(0)} \to M_i$ and $E_i^{(1)} \to M_i$ are isometrically isomorphic in a connection preserving way outside a compactum $K_i$,
		\item $\hat{A}(M_i,E_i) = \displaystyle\int_{M_i} \hat{A}(TM_i) \wedge (\mathrm{ch}(E_i^{(1)})-\mathrm{ch}(E_i^{(0)})) \neq 0$ for all $i \in \N$ and
		\item $\|R^{E_i}\|_\infty \lto 0$ for $i \lto \infty$.
	\end{itemize}
\end{Satz}
Roughly speaking, the construction is the following. When $\hat{A}(M) \neq 0$, then one can just take the constant sequence consisting of the identity $M \to M$ as covering and the trivial $\Z/2\Z$-graded bundle $\underline{\C} \oplus \underline{0} \to M$ as bundle. Else, if $\hat A(M)=0$, one may take a sequence $M_i \to M$ such that $M_i$ admits a $\frac{1}{i}$-area-contracting map $M_i \to S^{2\ell}$, where $2\ell = k > 0$ is chosen as in the definition of enlargeability. The bundles are obtained by pulling back a bundle $E^{(0)} \oplus E^{(1)} \to S^{2\ell}$, where $E^{(0)} \to S^{2\ell}$ satisfies $c_\ell(E^{(0)}) \neq 0$ and  $E^{(1)} \to S^{2\ell}$ is a trivial bundle of the same rank.

The second part consists of calculating the index of the Dirac operator on the twisted spinor bundle $\Sigma M_i \otimes E_i$ in two different ways. On the one hand, by the  the relative index theorem, its index is $\hat A(M_i,E_i) \neq 0$. On the other hand, assuming that $M$ carries a psc metric, the twisted Schrödinger-Lichnerowicz formula implies that this Dirac operator is invertible and thus has index zero, for large $i \in \N$. As we are interested in initial data sets, we replace this second step and obtain our main theorem:
\begin{Satz}[Main Theorem] \label{Thm:EnlObst}
	Let $M$ be a compact spin manifold that is $\hat{A}$-area-enlargeable.
	Then the path-components $C^-$ and $C^+$ of $\DEC(M)$ do not agree, i.\,e.\ if $g$ is a metric on $M$ and $\tau > 0$ is chosen so large that $(g,-\tau g), (g, \tau g) \in \DEC(M)$, then $(g,-\tau g)$ and  $(g, \tau g)$ belong to different path-components of $\DEC(M)$.
	\begin{proof}
		We first consider the case where the dimension of $M$ is even.
		Let $g$ be a metric on $M$ and $\tau > 0$ be large enough that $(g,-\tau g), (g, \tau g) \in \DEC(M)$.
		We choose a sequence of complex vector bundles $E_i \to M_i$ as in \cref{Thm:GL} and denote by $g_i$ the pull-back metric of $g$ on $M_i$. From the relative index theorem \cref{Cor:RelIndThm}, we obtain that for all $i \in \N$ the twisted index difference $\idiff^{E_i}((g_i, -\tau g_i),(g_i, \tau g_i))$ corresponds to $\hat A(M_i, E_i) \neq 0$ under the isomorphism $\mathrm{K}^1(I, \del I) \cong \Z$, in particular it is non-zero.
		
		We now assume for contradiction that  $(g,-\tau g)$ and $(g, \tau g)$ are connected in $\DEC(M)$ by a path $t \mapsto (g(t),k(t))$. As the interval $I$ is compact, there is a constant $c > 0$, such that $\rho(t)-\|j(t)\| \geq 4c$ for all $t \in I$. Of course, this holds as well for the pulled-back path in $\DEC(M_i)$, with the same constant. Since for any $\phi \in \overline{\Sigma} M_i$ and $e \in E_i$
		\begin{align*}
			\|\mathcal{R}^{E_i} (\phi \otimes e)\| \leq \sum_{j<k} \|\phi\| \|R^{E_i}(e_j,e_k) e\| \leq \frac{n(n-1)}{2} \|R^{E_i}\| \|\phi \otimes e\|.
		\end{align*}
		and $\|R^{E_i}\|_\infty \lto 0$ for $i \lto \infty$, we have $\|\mathcal{R}^{E_i}\| < c$ as long as $i \in \N$ is large enough. Hence for large $i \in \N$ the pulled back path $t \mapsto (g_i(t), k_i(t))$ lies entirely in $\Ini^{> \mathcal{R}^{E_i} + c}(M_i)$. Thus by \cref{Prop:DefInd} $\idiff^{E_i}((g_i, -\tau g_i),(g_i, \tau g_i)) = 0$, which is the desired contradiction.

		In the odd-dimensional case, we replace $M$ by $M \times S^1$, which will again be $\hat A$-area-enlargeable and spin, and is of even dimension. Thus, we conclude that $C^+ \neq C^-$ in $\DEC(M \times S^1)$. By \cref{Lem:Stab}, the same holds for $\DEC(M)$.
	\end{proof}
\end{Satz}

\appendix
\section{Schrödinger-Lichnerowicz formula for the twisted Dirac-Witten operator}
\begin{Satz}
	For all $\psi \in \Gamma((\Sigma M \oplus \Sigma M) \otimes_\C E)$
	\begin{align*}
	\left(\overline{D}^E\right)^2 \psi = \overline{\nabla}^* \overline{\nabla} \psi + \frac{1}{2} (\rho - e_0 \cdot j^\sharp \cdot) \psi + \mathcal{R}^E \psi,
	\end{align*}
	where $\rho$ and $j$ are defined as in \eqref{eq:CE} in terms of the pair $(g,k)$ and $\mathcal{R}^E(\phi \otimes e)= \sum_{i<j} e_i \cdot  e_j \cdot \phi \otimes R^E(e_i,e_j)e$ for $\phi \otimes e \in (\Sigma_p M \oplus \Sigma_p M) \otimes_\C E_p$ and an orthonormal basis $(e_1, \ldots, e_n)$ of $T_pM$, $p \in M$.
\begin{proof}
	We show how to reduce the formula to the Schödinger-Lichnerowicz type formula in the untwisted case \eqref{eq:SL}, using a local calculation. For this, let $(e_1,\ldots, e_n)$ be a local orthonormal frame. Without loss of generality, we may assume that $\psi$ can be written locally as $\phi \otimes e$ as everything is linear. Then
	\begin{align*}
		\left(\overline{D}^E\right)^2(\phi \otimes e) &= 
		\sum_{i,j} e_i \cdot \overline{\nabla}_{e_i}(e_j \cdot \overline{\nabla}_{e_j} \phi) \otimes e 
		+ \sum_{i,j} e_i \cdot \overline{\nabla}_{e_i}(e_j \cdot \phi) \otimes \nabla_{e_j}^E e \\
		&\phantom{=}\;+ \sum_{i,j} e_i \cdot e_j \cdot (\overline{\nabla}_{e_j} \phi) \otimes \nabla_{e_i}^E e 
		+ \sum_{i,j} e_i \cdot e_j \cdot \phi \otimes \nabla_{e_i}^E \nabla_{e_j}^E e \\
		&= (\overline{D}^2 \phi) \otimes e 
		+ \sum_{i,j} e_i \cdot (\overline{\nabla}_{e_i} e_j) \cdot \phi \otimes \nabla_{e_j}^E e \\
		&\phantom{=}\,-2 \sum_{i} (\overline{\nabla}_{e_i} \phi) \otimes \nabla_{e_i}^E e 
		+ \sum_{i,j} e_i \cdot e_j \cdot \phi \otimes \nabla_{e_i}^E \nabla_{e_j}^E e
	\intertext{and}
		\overline{\nabla}^*\overline{\nabla}(\phi \otimes e) &=
		\sum_{i} \overline{\nabla}^*(e_i^* \otimes \overline{\nabla}_{e_i} (\phi \otimes e)) \\
		&= -\sum_{i} \overline{\nabla}_{e_i} \overline{\nabla}_{e_i} (\phi \otimes e) -\sum_{i} e_0 \cdot k(e_i,-)^\sharp \cdot \overline{\nabla}_{e_i} (\phi \otimes e) + \sum_{i} \overline{\nabla}_{\nabla_{e_i} e_i} (\phi \otimes e) \\
		&= -\sum_{i} (\overline{\nabla}_{e_i} \overline{\nabla}_{e_i} \phi) \otimes e - 2 \sum_{i} (\overline{\nabla}_{e_i} \phi) \otimes \nabla_{e_i}^E e - \sum_{i} \phi \otimes \nabla_{e_i}^E \nabla_{e_i}^E e \\
		&\phantom{=}\; -\sum_{i} e_0 \cdot k(e_i,-)^\sharp \cdot (\overline{\nabla}_{e_i} \phi) \otimes e -\sum_{i} e_0 \cdot k(e_i,-)^\sharp \cdot  \phi \otimes \nabla_{e_i}^E e  \\
		&\phantom{=}\; + \sum_{i} (\overline{\nabla}_{\nabla_{e_i} e_i} \phi) \otimes e + \sum_{i} \phi \otimes \nabla_{\nabla_{e_i} e_i}^E e \\
		&= (\overline{\nabla}^* \overline{\nabla} \phi) \otimes e - 2 \sum_{i} (\overline{\nabla}_{e_i} \phi) \otimes \nabla_{e_i}^E e - \sum_{i} \phi \otimes \nabla_{e_i}^E \nabla_{e_i}^E e  \\
		&\phantom{=}\;-\sum_{i} e_0 \cdot k(e_i,-)^\sharp \cdot  \phi \otimes \nabla_{e_i}^E e  + \sum_{i} \phi \otimes \nabla_{\nabla_{e_i} e_i}^E e	
	\end{align*}
	using $\overline{\nabla}_X \psi = \nabla_X \psi - \frac{1}{2} e_0 \cdot k(X,-)^\sharp \cdot \psi$ and that the formal adjoint of $\nabla$ is given by $\nabla^* \colon \alpha \otimes \psi \mapsto - \sum_j \nabla_{e_j}(\alpha \otimes \psi)(e_j) = -\sum_j \nabla_{e_j} (\alpha(e_j)\psi) +\sum_j \alpha(\nabla_{e_j} e_j) \psi$, $\alpha \in \Omega^1(M)$.
	Noting that
	\begin{align*}
		\sum_{i,j} e_i \cdot (\overline{\nabla}_{e_i} e_j) \cdot \phi \otimes \nabla_{e_j}^E e &= \sum_{i,j,k}  g(\nabla_{e_i} e_j, e_k) e_i \cdot e_k \cdot \phi \otimes \nabla_{e_j}^E e + \sum_{i,j} k(e_i,e_j) e_i \cdot e_0 \cdot \phi \otimes \nabla_{e_j}^E e  \\
		&= -\sum_{i,j,k} g(e_j, \nabla_{e_i} e_k) e_i \cdot e_k \cdot \phi \otimes \nabla_{e_j}^E e + \sum_{j} k(-,e_j)^\sharp \cdot e_0 \cdot \phi \otimes \nabla_{e_j}^E e  \\
		&= -\sum_{i,j} e_i \cdot e_j \cdot \phi \otimes \nabla_{\nabla_{e_i} e_j}^E e - \sum_{i} e_0 \cdot k(e_i,-)^\sharp \cdot \phi \otimes \nabla_{e_i}^E e
	\end{align*}
	this implies
	\begin{align*}
		\left(\overline{D}^E\right)^2(\phi \otimes e) - \overline{\nabla}^*\overline{\nabla}(\phi \otimes e) &=
		(\overline{D}^2 \phi) \otimes e  - (\overline{\nabla}^* \overline{\nabla} \phi) \otimes e -\sum_{i \neq j} e_i \cdot e_j \cdot \phi \otimes \nabla_{\nabla_{e_i} e_j}^E e\\
		&\phantom{=}\; + \sum_{i \neq j} e_i \cdot e_j \cdot \phi \otimes \nabla_{e_i}^E \nabla_{e_j}^E e.
	\end{align*}
	Using the untwisted Schrödinger-Lichnerowicz type formula \eqref{eq:SL}, the first two terms compute to $\frac{1}{2}(\rho-e_0 \cdot j^\sharp \cdot ) \phi \otimes e$. Thus it remains to identify the remaining terms with $\mathcal{R}^E(\phi \otimes e)$:
	\begin{align*}
		\mathcal{R}^E(\phi \otimes e) &=  \sum_{i<j} e_i \cdot  e_j \cdot \phi \otimes R^E(e_i,e_j)e \\
			&= \sum_{i<j} e_i \cdot  e_j \cdot \phi \otimes (\nabla_{e_i}^E \nabla_{e_j}^E - \nabla_{e_j}^E \nabla_{e_i}^E - \nabla_{\nabla_{e_i} e_j}^E + \nabla_{\nabla_{e_j}e_i}^E)e \\
			&= \sum_{i \neq j} e_i \cdot  e_j \cdot \phi \otimes (\nabla_{e_i}^E \nabla_{e_j}^E - \nabla_{\nabla_{e_i} e_j}^E)e. \qedhere
	\end{align*}
\end{proof}
\end{Satz}

\addcontentsline{toc}{section}{References}


\end{document}